\documentclass{amsart}

\usepackage{amsfonts}
\usepackage{amssymb}
\usepackage{latexsym}

\newtheorem{theorem}{Theorem}[section]
\newtheorem{lemma}[theorem]{Lemma}

\theoremstyle{definition}

\theoremstyle{remark}

\numberwithin{equation}{section}

\newcommand{\Q}{\mathbb{Q}}
\newcommand{\R}{\mathbb{R}}
\newcommand{\Z}{\mathbb{Z}}

\begin{document}

\title[Nonvanishing Vector Fields on Almost-Complex, Cyclic Orbifolds]
{A Complete Obstruction to the Existence of
    Nonvanishing Vector Fields on Almost-Complex, Closed, Cyclic Orbifolds}

\author{Christopher Seaton}
\address{Department of Mathematics and Computer Science,
Rhodes College, 2000 N. Parkway,
Memphis, TN 38112}
\email{seatonc@rhodes.edu}


\date{August 2nd, 2006}


\begin{abstract}
We determine several necessary and sufficient conditions for a
closed almost-complex orbifold $Q$ with cyclic local groups to admit
a nonvanishing vector field. These conditions are stated separately
in terms of the orbifold Euler-Satake characteristics of $Q$ and its
sectors, the Euler characteristics of the underlying topological
spaces of $Q$ and its sectors, and in terms of the orbifold Euler
class $e_{orb}(Q)$ in Chen-Ruan orbifold cohomology $H_{orb}^\ast
(Q; \R)$.
\end{abstract}

\maketitle


\section{Introduction}

Orbifolds are singular spaces locally modeled by $\R^n / G$ where
$G$ is a finite subgroup of $O(n)$ that acts with a fixed-point set
of codimension at least 2. The original definition of an orbifold
was introduced by Satake in \cite{satake1} under the name
$V$-manifold, and the term orbifold was given by Thurston in
\cite{thurston}. Thurston's orbifolds included a larger class than
those of Satake, for he allows the local groups $G$ to act with a
fixed-point set of codimension 1.  Today, the definition of an
orbifold varies from author to author. Here, we retain the
requirement that the local groups act with a fixed-point set of
codimension at least 2, but do not require the groups to act
effectively.  Hence, Satake's $V$-manifolds correspond to our {\bf
reduced} orbifolds.

Let $Q$ be a closed, reduced orbifold of dimension $n$.  One of the first
things that was studied on orbifolds is the generalization of de Rham theory by
Satake in \cite{satake1} and \cite{satake2}.
In the latter of these two papers, Satake developed a generalization
of the Poincar\'{e}-Hopf Theorem, that if $X$ is a vector field on $Q$ with only
isolated zeros, then
\begin{equation}
\label{eq-ph}
    \mbox{ind}_{orb}(X) = \chi_{orb}(Q).
\end{equation}
Here, $\mbox{ind}_{orb}(X)$ is the orbifold index of the vector
field and $\chi_{orb}(Q)$ the orbifold Euler-Satake characteristic
of $Q$ (see Section \ref{sec-defs} for a review of the definitions).

More recently, the author has developed an additional generalization of the
Poincar\'{e}-Hopf Theorem to orbifolds \cite{mythesis}.  In this case, the left
side of  the equation is the orbifold index of the vector field $\tilde{X}$ induced
by $X$ on $\tilde{Q}$, the space of sectors of the orbifold.  The right side
then becomes $\chi(\mathbb{X}_Q)$, the Euler characteristic of the underlying
topological space $\mathbb{X}_Q$ of $Q$:
\begin{equation}
\label{eq-myph}
    \mbox{ind}_{orb}(\tilde{X}) = \chi(\mathbb{X}_Q).
\end{equation}
We will review these definitions in the sequel;
here, we note that if $X$ is a nonvanishing vector field, then $\tilde{X}$ is
nonvanishing as well.

As in the case of manifolds \cite{phopf}, it is a direct corollary
of these formulae that an orbifold admits a nonvanishing vector
field only if its orbifold Euler-Satake characteristic vanishes (in
the case of Equation \ref{eq-ph}), and the Euler characteristic of
its underlying topological space vanishes (in the case of Equation
\ref{eq-myph}).  Unlike the case of manifolds, however, the converse
of both of these statements is false.  It is easy to construct
examples of $2$-orbifolds $Q$ such that $\chi_{orb}(Q) = 0$ or
$\chi( \mathbb{X}_Q) = 0$, yet whose singular points force any
vector field to vanish.  While it is impossible for both of these
invariants to vanish for a nontrivial $2$-orbifold, it is possible
to construct a $4$-dimensional orbifold such that $\chi_{orb}(Q) =
\chi( \mathbb{X}_Q) = 0$ that does not admit a nonvanishing vector
field. For instance, one may take an orbifold whose underlying space
is $\mathbb{T}^4$ and whose singular set is the disjoint union of
$S^2$ and a surface of genus $2$, all with isotropy group $\Z_3$.

In this paper, we determine necessary and sufficient conditions for
a closed, almost-complex orbifold with cyclic local groups to admit
a nonvanishing vector field. Our main result is the following
theorem.


\begin{theorem}
\label{thrm-mainresult} Let $Q$ be a closed almost-complex cyclic
orbifold, and then the following are equivalent:

\vspace{0.1cm}

(i) $Q$ admits a nonvanishing vector field.

(ii)    $\tilde{Q}$ admits a nonvanishing vector field.

(iii)   The Euler characteristic of the underlying space of each
    sector $\tilde{Q}_{(g)}$ is zero.

(iv)    The orbifold Euler-Satake characteristic
    of each sector $\tilde{Q}_{(g)}$ is zero.

(v) $e_{orb}(Q)$, the orbifold Euler class of $Q$, is zero in
$H_{orb}^\ast(Q ; \R)$.

\end{theorem}

In Section \ref{sec-defs}, we review the pertinent definitions and
fix our notation.  The main constructions we require are that of the
space of sectors of an orbifold, Chen-Ruan orbifold cohomology, and
the orbifold Euler class; the reader is referred to the original
sources for a more detailed exposition.  In Section
\ref{sec-structure}, we study the relationship between the sectors
of an orbifold.  Section \ref{sec-mainresult} contains the proof of
our theorem.

The author is pleased to acknowledge Carla Farsi, Alexander Gorokhovsky,
Judith Packer, Arlan Ramsay, and Lynne Walling
for useful discussions and support during the work leading to this result.


\section{Review of Definitions}
\label{sec-defs}

In this section, we briefly review the definitions we will need.
For more information,
the reader is referred to the original work of Satake in
\cite{satake1} and \cite{satake2}.  As well,
\cite{ruangwt} contains as an appendix a thorough introduction to orbifolds,
focusing on their differential geometry, and \cite{mythesis} contains
an introduction to orbifolds with an emphasis on vector fields.

A (${\mathcal C}^\infty$) {\bf orbifold} $Q$ is a Hausdorff space
$\mathbb{X}_Q$ such that each point is contained
in an open set modeled by an {\bf orbifold chart} or {\bf local
uniformizing system.}  By this, we mean a triple $\{ V, G, \pi \}$ where
\begin{itemize}
\item       $V$ is an open subset of $\R^n$,
\item       $G$ is a finite group with a ${\mathcal C}^\infty$
        action on $V$ such that the fixed point
        set of any $\gamma \in G$ which does not
        act trivially on $V$ has codimension at least 2 in $V$, and
\item       $\pi : V \rightarrow U$ is a surjective continuous map such
        that $\forall \gamma \in G$, $\pi \circ \gamma = \pi$ that
        induces a homeomorphism $\tilde{\pi} : V/G \rightarrow U$.
\end{itemize}
The image $U=\pi(V)$ is called a {\bf uniformized set} in $Q$. The
group $G$ is known as a {\bf local group.}  If the local group of a
chart $\{ V, G, \pi \}$ acts effectively, then the chart is said to
be {\bf reduced}; if all charts are reduced, then $Q$ is a {\bf
reduced orbifold.}  In the spirit of \cite{chenhu} for the case of
Abelian orbifolds, we adopt the convention that if each local group
is cyclic, then $Q$ is a {\bf cyclic orbifold}.

It is required that if a point $p$ is contained in two uniformized
sets $U_i$ and $U_j$, then there is a uniformized set $U_k$ such
that $p \in U_k \subset U_i \cap U_j$.  Moreover, if $U_i \subseteq
U_j$ are two sets uniformized by $ \{ V_i , G_i , \pi_i \}$ and $ \{
V_j , G_j , \pi_j \}$, respectively, then we require that they are
related by an injection $\lambda_{ij} : \{ V_i , G_i , \pi_i \}
\rightarrow \{ V_j , G_j , \pi_j \}$.  An {\bf injection}
$\lambda_{ij}$ is a pair $\{ f_{ij}, \phi_{ij} \}$ where
\begin{itemize}
\item       $f_{ij} : G_i \rightarrow G_j$ is an injective
homomorphism such
        that if $K_i$ and $K_j$ denote the kernel of the action of $G_i$ and $G_j$,
        respectively, then $f_{ij}$ restricts to an isomorphism of $K_i$ onto $K_j$,
        and
\item       $\phi_{ij} : V_i \rightarrow V_j$ is a smooth embedding
such that
        $\pi_i = \pi_j \circ \phi_{ij}$ and such that for each $\gamma \in G_i$,
        $\phi_{ij} \circ \gamma = f_{ij}(\gamma) \circ \phi_{ij}$
\end{itemize}
(see \cite{satake2}, \cite{ruangwt}, or \cite{mythesis}).

Orbifold vector bundles are defined over each chart $\{ V, G, \pi
\}$ as $G$-vector bundles over $V$.  In particular, the tangent
bundle is defined locally to be the ordinary tangent bundles $TV$
with $G$-structure given by the differential of the $G$-action.
Sections of orbifold vector bundles correspond locally to
$G$-equivariant sections of the $G$-bundles over $V$.

In particular, with the help of the exponential map, it is possible to
replace each orbifold chart containing a point $p \in Q$ with an equivalent
chart such that $G$ acts on $V$ as a subgroup of $O(n)$, and $p$
is the image under $\pi$ of the origin in $V$.  Such a chart will
be known as a {\bf chart at $p$}, denoted $\{ V_p, G_p, \pi_p \}$,
with $U_p := \pi_p(V_p)$, etc.  Note that in a chart at $p$,
$G_p$ is the isotropy group of $p$, and its isomorphism class is
independent of the choice of chart.

The {\bf orbifold index of a vector field $X$ on $Q$} with isolated
zeros is defined to be the sum of the indices at each zero of $X$,
where the index at a zero point $p$ is the quotient of the (usual)
index of the vector field in an orbifold chart and the order of the
isotropy group at $p$.  In other words, if $\{ V_p, G_p, \pi_p \}$
is an orbifold chart at $p$, then the index of $X$ at $p$ is
$\frac{1}{|G_p|} \mbox{ind}_{\mathbf{0}} (\pi_p^\ast X)$.

The {\bf orbifold Euler-Satake characteristic} $\chi_{orb}(Q)$ is most easily defined
by finding a simplicial decomposition ${\mathcal T}$ for $Q$ such that the
isomorphism class of the isotropy group of each point on the interior of a simplex
is constant (such a triangulation always exists; see \cite{moerdijk}).
For each simplex  $\sigma \in {\mathcal T}$, if we let $m_{\sigma}$ denote the order of
this isotropy group, then
\[
    \chi_{orb}(Q) := \sum\limits_{\sigma \in {\mathcal T}}
        (-1)^{\mbox{dim}\:\sigma} \frac{1}{m_{\sigma}}.
\]
Note that if $A \cup B$ is a union of orbifolds, then it is
straightforward to show from this definition that, as in the case of
the usual Euler characteristic,
\begin{equation}
\label{eq-addativeesk} \chi_{orb}(A \cup B) = \chi_{orb}(A) +
\chi_{orb}(B) - \chi_{orb}(A \cap B).
\end{equation}

The Chen-Ruan orbifold cohomology groups are defined in terms of the
{\bf space of sectors} of the orbifold (also known as the {\bf
inertia orbifold}). We recall the construction, referring the reader
to \cite{chenruan} and \cite{kawasaki2} for more details.

Select for each $p \in Q$ a chart $\{ V_p, G_p, \pi_p \}$ at $p$.
Then the set
\[
    \tilde{Q} := \{ (p, (g)_{G_p} ) : p \in Q, g \in G_p \}
\]
(where $(g)_{G_p}$ is the conjugacy class of $g$ in $G_p$) is naturally
an orbifold, with local charts
\[
    \{ \pi_{p, g} : (V_p^g, C(g) ) \rightarrow V_p^g / C(g)
        : p \in Q, g \in G_p \} .
\]
Here, $V_p^g$ is the fixed point set of $g$ in $V_p$, and $C(g)$ is
the centralizer of $g$ in $G_p$.

An equivalence relation can be placed on the conjugacy classes of
the local groups $G_p$ as follows.  The conjugacy class of $g \in
G_p$ and that of $h \in G_q$ are equivalent if there is an injection
of a chart at $q$ into a chart at $p$ such that $f_{qp}(h) = g$ for
the corresponding homomorphism $f_{qp}$.  Note that, as the choice
of injection is not generally unique, this equivalence is defined on
conjugacy classes. In the case that $G_p$ and $G_q$ are Abelian, of
course, each conjugacy class contains one element (see \cite{chenhu}
for more details in the case of Abelian local groups).

Let $T$ denote the set of equivalence classes under this relation
(which is finite for $Q$ compact) and $(g)$ the equivalence class of
a conjugacy class $(g)_{G_p}$.  Then
\[
    \tilde{Q} = \bigsqcup\limits_{(g) \in T } \tilde{Q}_{(g)} ,
\]
where
\[
    \tilde{Q}_{(g)} = \{ (p, (g^\prime)_{G_p})
        : g^\prime \in G_p , (g^\prime)_{G_p} \in (g) \} .
\]

Each of the $\tilde{Q}_{(g)}$ for $(g) \neq (1)$ is called a {\bf
twisted sector}; $\tilde{Q}_{(1)}$ is the {\bf nontwisted sector},
and is diffeomorphic to $Q$ as an orbifold.  The {\bf sectors} of
the orbifold refer to both the twisted sectors and the nontwisted
sector.  We note that even in the case that $Q$ is connected, a
twisted sector of $Q$ need not be.

If $Q$ is an almost complex orbifold, a function $\iota : \tilde{Q}
\rightarrow \Q$ is defined which is constant on each $\tilde{Q}$.
The value of this function on $\tilde{Q}_{(g)}$, denoted
$\iota_{(g)}$, is called the degree shifting number of $(g)$. The
orbifold cohomology groups are defined by
\[
    H_{orb}^d (Q; \R)   :=  \bigoplus_{(g) \in T}
        H^{d - 2\iota_{(g)}} (\tilde{Q}_{(g)} ; \R),
\]
where the groups on the right side are the usual de Rham cohomology
groups of the orbifolds $\tilde{Q}_{(g)}$.

For a vector bundle $\rho :E \rightarrow Q$, the space $E$ is
naturally an orbifold, so it is possible to form $\tilde{E}$ as
above.  In the case that $E$ is a {\bf good} vector bundle (see
\cite{ruangwt}, Section 4.3), $\tilde{E}$ is naturally an orbifold
vector bundle over $\tilde{Q}$, although its dimension varies in
general over the connected components of $\tilde{Q}$.  Similarly,
smooth sections $s : Q \rightarrow E$ of a good vector bundle $E$
naturally induce smooth sections $\tilde{s} : \tilde{Q} \rightarrow
\tilde{E}$ in such a way that if $s$ is nonvanishing, then so is
$\tilde{s}$ (see \cite{mythesis}, Lemma 4.4.1; the definition of an
orbifold vector bundle in this reference corresponds to Ruan's
definition of a good orbifold vector bundle). In particular, this is
true of the orbifold tangent bundle $TQ$, which is always a good
vector bundle. Moreover, $T\tilde{Q} = \widetilde{TQ}$, so that
nonvanishing vector fields over $Q$ naturally induce nonvanishing
vector fields over $\tilde{Q}$.

A connection on a good vector bundle $E$ induces one on $\tilde{E}$.
The {\bf orbifold Euler class} of an orbifold vector bundle $E$ is
defined in terms of such a connection.  It can be taken to be the
sum in $H^\ast_{orb}(Q ; \R)$ of the usual Euler classes of the
orbifold bundle $\tilde{E}$ restricted to each of the connected
components of each $\tilde{Q}_{(g)}$ (see \cite{mythesis}).

We will require the formula
\begin{equation}
\label{eq-2ndgb}
    \sum\limits_{(g) \in T} \chi_{orb} ( \tilde{Q}_{(g)} )
    =
    \chi (\mathbb{X}_Q).
\end{equation}
See the proof of the Second Gauss-Bonnet Theorem for orbifolds
(Theorem 4.4.2) in \cite{mythesis} for a verification of this
formula.


\section{The Structure of the Sectors of a Cyclic Orbifold}
\label{sec-structure}

The construction of $\tilde{Q}$ decomposes an orbifold $Q$ into
multiple sectors, the largest of these being diffeomorphic to $Q$
and the others being simpler orbifolds of lower dimension.  In this
section, we study the structure of these connected components as
they appear as subsets of a cyclic orbifold $Q$.

Let $\pi : \tilde{Q} \rightarrow Q$ denote the projection with
$\pi(p, (g)) = p$.  Then $\pi$ is a ${\mathcal C}^\infty$ map (see
\cite{chenruan}). In the case that each of the local groups is
Abelian, clearly $C(g) = G_p$ for each $(g) \in T$, so that the
local uniformized sets $V_p^g / C(g)$ are diffeomorphic to subsets
of uniformized sets $V_p/G$ in $Q$ via $\pi$. Therefore, in this
case, $\pi$ is an embedding when restricted to $\tilde{Q}_{(g)}$ for
each $(g) \in T$ (see \cite{chenhu}).  We define the following
relation on the $\tilde{Q}_{(g)}$ via $\pi$: we say that
$\tilde{Q}_{(h)} \leq \tilde{Q}_{(g)}$ whenever
$\pi(\tilde{Q}_{(h)}) \subseteq \pi(\tilde{Q}_{(g)})$,
$\tilde{Q}_{(h)} \equiv \tilde{Q}_{(g)}$ whenever
$\pi(\tilde{Q}_{(h)}) = \pi(\tilde{Q}_{(g)})$, etc.  Note that
$\tilde{Q}_{(h)} \equiv \tilde{Q}_{(g)}$ does not imply that
$\tilde{Q}_{(h)} = \tilde{Q}_{(g)}$, but that, through the
appropriate restrictions of $\pi$ and inverses of these restrictions
(see Equation \ref{eq-restrictdiffeo}), they are diffeomorphic.
Therefore, this relation can be thought of as a partial order of
equivalence classes of the sectors under the equivalence relation
$\equiv$. We will refer to the elements of the minimal equivalence
classes under $\equiv$ as minimal with respect to the relation
$\leq$. As well, it will be convenient for us to state this relation
in terms of the elements of $T$; i.e. $(h) \leq (g)$ will mean that
$\tilde{Q}_{(h)} \leq \tilde{Q}_{(g)}$, etc.

With respect to this relation, the nontwisted sector
$\tilde{Q}_{(1)}$ is clearly maximal, as are each of the
$\tilde{Q}_{(g)}$ where the representatives of $(g)$ act trivially
(in the case that $Q$ is not reduced).  Similarly, we have the
following:

\begin{lemma}
\label{lemm-minmnflds}

Let $Q$ be a closed, cyclic orbifold, and let $\tilde{Q}_{(g)}$ be a
sector of $Q$ that is minimal with respect to the relation given
above.  Then $\tilde{Q}_{(g)}$ is a manifold equipped with the
trivial action of a finite group, and hence the associated reduced
orbifold $(\tilde{Q}_{(g)})_{red}$ is a manifold.

\end{lemma}

\begin{proof}

Let $\tilde{Q}_{(g)}$ be a minimal sector, and let $(p, (g))$ be a
point in $\tilde{Q}_{(g)}$.  Fix an orbifold chart at $(p, (g))$ of
the form $\{ V_p^g, C(g), \pi_{p, g} \}$ induced by an orbifold
chart $\{ V_p, G_p, \pi_p \}$ for $Q$ at $p$, and note that as $G_p$
is Abelian, $C(g) = G_p$.  Without loss of generality, assume that
the representative $g$ of $(g)$ is an element of $G_p$.

The subgroup $\langle g \rangle$ of $G_p$ generated by $g$ clearly
acts trivially on $V_p^g$, and a representative of $(g)$ is in the
isotropy group of each point in $\pi (\tilde{Q}_{(g)})$.  We argue
that $\langle g \rangle$ is exactly the isotropy group of each point
in $\pi (\tilde{Q}_{(g)})$.

Let $(q, (g_q))$ be an arbitrary point in $\tilde{Q}_{(g)}$, and
then $g_q$ is a representative of $(g)$ in the isotropy group $G_q$.
Hence, there is a sequence of points $(p_i, (g_{p_i}))$, $i = 0,
\ldots , s$ such that $p_0 = p$, $p_s = q$, $p_i$ and $p_{i+1}$ lie
in the same chart $\{ V_i, G_i, \pi_i \}$, and $g_i = g_{i+1}$ in
$G_i$ (see \cite{chenhu}; note that in general, $g_i$ and $g_{i+1}$
are required to be conjugate in $G_i$, but that this implies that
they are equal in our case).  Fix a chart $\{ V_q, G_q, \pi_q \}$ at
$q \in Q$, and let $a$ be a generator of $G_q$.  Then there is a $k
\in \Z$ such that $a^k = g_q$.

It is obvious that $V_q^a \subseteq V_q^{g_q}$.  Moreover, for any
injection of $\{ V_q, G_q, \pi_q \}$ (or a restriction to a smaller
chart containing $q$) into another chart, the associated injective
homomorphism $f$ satisfies $f(a)^k = f(g_q)$. This implies that in
any chart for $Q$, a representative of $(a)$ has a fixed-point set
contained in that of a representative of $(g_q)$.  Therefore,
$\pi(\tilde{Q}_{(a)})$ is contained in $\pi(\tilde{Q}_{(g_q)})$.  As
$(g_q)$ and $(g)$ represent the same equivalence class in $T$,
$\tilde{Q}_{(g_q)} = \tilde{Q}_{(g)}$. This implies that $(a) \leq
(g)$.  Hence, as $\tilde{Q}_{(g)}$ is assumed to be minimal, we must
have that $(a) \equiv (g)$. Therefore, $a = g_q \in G_q$.

Now, for each pair $p_i, p_{i+1}$ of points above, the intersection
of the chart at $p_i$ with that at $p_{i+1}$ produces a chart at
$p_i$ contained in that at $p_{i+1}$.  Hence, there is an injection
of charts including an injective homomorphism $f_{i, i+1} : G_{p_i}
\rightarrow G_{p_{i+1}}$, producing a sequence
\[
    G_p = G_{p_0} \stackrel{f_{0, 1}}{\longrightarrow} G_{p_1} \stackrel{f_{1, 2}}{\longrightarrow}
        \cdots \stackrel{f_{s-1, s}}{\longrightarrow} G_{p_s} = G_q.
\]
Let $F: G_p \rightarrow G_q$ be the composition of these maps, and
then note that, as $a = g_q$ is a generator of $G_q$, and as $F(g) =
g_q$, we must have that each $f_{i, i+1}$ is an isomorphism.  In
particular, $g$ generates $G_p$, and $G_q$ is isomorphic to $G_p$.

As $(q, (g_q))$ was arbitrary, we conclude that the isotropy group
of each point in $\tilde{Q}_{(g)}$ is isomorphic to $G_p = \langle g
\rangle$.  Therefore, $\tilde{Q}_{(g)}$ is a manifold equipped with
the trivial action of the finite group $G_p$.

\end{proof}

In the case of a cyclic orbifold, the twisted sectors each decompose
into twisted sectors themselves.  In particular, as was noted above,
for each $(g) \in T$, the map
\[
    \pi_{(g)} := \pi_{|\tilde{Q}_{(g)}} : \tilde{Q}_{(g)} \rightarrow
    Q
\]
is an embedding of $\tilde{Q}_{(g)}$ into $Q$.  Hence, if $(h) \leq
(g)$, then as $\tilde{Q}_{(h)}$ is embedded into
$\pi(\tilde{Q}_{(g)})$, the composition
\begin{equation}
\label{eq-restrictdiffeo}
    \tilde{Q}_{(h)} \stackrel{\pi_{(h)}}{\rightarrow}
    \pi(\tilde{Q}_{(h)}) \stackrel{\pi_{(g)}^{-1}}{\rightarrow}
    \tilde{Q}_{(g)}
\end{equation}
defines an embedding of $\tilde{Q}_{(h)}$ into $\tilde{Q}_{(g)}$ (we
note that the inverse $\pi_{(g)}^{-1}$ is defined on
$\pi(\tilde{Q}_{(g)}) \supseteq \pi_{(h)}(\tilde{Q}_{(h)})$, and not
on $Q$). Denote this composition
\[
    \pi_{(h), (g)} := \pi_{(g)}^{-1} \circ \pi_{(h)} : \tilde{Q}_{(h)}
    \rightarrow \tilde{Q}_{(g)}.
\]
If $(g) \equiv (h)$, then $\pi_{(h), (g)}$ is a diffeomorphism of
orbifolds.

We conclude this section with the following Lemma, which illustrates
that the singular set of a sector decomposes into the image of
embeddings of strictly smaller sectors.  Note that we reduce the
sector $\tilde{Q}_{(g)}$ only to distinguish between the singular
and regular points of the reduced orbifold; for $(g) \neq (1)$, all
of $\tilde{Q}_{(g)}$ is singular.

\begin{lemma}
\label{lemm-singinsector}

Let $Q$ be a closed, cyclic orbifold.  Fix a sector
$\tilde{Q}_{(g)}$, and let $(p, (g))$ be a singular point of
$(\tilde{Q}_{(g)})_{red}$. Then there is a $(g^\prime) \in T$ with
$(g^\prime) < (g)$ such that $(g^\prime)$ has a representative in
$G_p$. Hence, each singular point of $(\tilde{Q}_{(g)})_{red}$ is
contained in the image of some such embedding $\pi_{(g^\prime),
(g)}$ of such a $\tilde{Q}_{(g^\prime)}$ into $\tilde{Q}_{(g)}$.
Moreover, the intersection of two such embeddings is the image of a
sector.

\end{lemma}

\begin{proof}

Fix a chart $\{ V_p, G_p, \pi_p \}$ for $Q$ at $p$, and then let $\{
V_p^g, C(g), \pi_{p, g} \}$ be the induced chart for
$\tilde{Q}_{(g)}$.  Assume without loss of generality that the
representative $g$ of $(g)$ is an element of $G_p$.  As $(p, (g))$
is assumed to be singular, there is an $h \in C(g) = G_p$ that acts
nontrivially on $V_p^g$ and fixes $p$. Hence, $h$ has a fixed point
subset $V_p^h$ in $V_p$ with $V_p^h \cap V_p^g \subsetneq V_p^g$.

As $G_p$ is cyclic, the subgroup $\langle g, h \rangle$ of $G_p$ is
cyclic as well.  Let $g^\prime$ be a generator of $\langle g, h
\rangle$, and then it is obvious that $V_p^{g^\prime} = V_p^g \cap
V_p^h$.

With this, we need only note that for any injection of a chart at
$q$ to another chart, the corresponding injective homomorphism $f$
satisfies $\langle f(g^\prime) \rangle = \langle f(g), f(h)
\rangle$, implying that the above relationships hold for any chart
containing representatives of $(g)$, $(h)$, and $(g^\prime)$.
Therefore, it follows that $(g^\prime) \leq (g)$ and $(g^\prime)
\leq (h)$.  In fact, we see that $\pi(\tilde{Q}_{(g)}) \cap
\pi(\tilde{Q}_{(h)}) = \pi(\tilde{Q}_{(g^\prime)})$, and as $h$ acts
nontrivially on $V^g$, we conclude that $(g^\prime) < (g)$.
Moreover, as $(g^\prime)$ has a representative in the isotropy group
$G_p$ of $p$, the singular point $(p, (g))$ lies in the image of
$\tilde{Q}_{(g^\prime)}$ under $\pi_{(g^\prime),(g)}$.

Now, suppose that $\tilde{Q}_{(h_1)}$ and $\tilde{Q}_{(h_2)}$ are
two sectors with $(h_1) < (g)$ and $(h_2) < (g)$. If the embeddings
of $\tilde{Q}_{(h_1)}$ and $\tilde{Q}_{(h_2)}$ intersect in
$\tilde{Q}_{(g)}$, then there is a point $(p, (g))$ such that $G_p$
contains representatives of $(h_1)$ and $(h_2)$ (say $h_1$ and
$h_2$). Again, $\langle g, h_1, h_2 \rangle$ is cyclic, so let
$g^\prime$ be a generator.  Repeating the above argument, it is
clear that
\[
\pi(\tilde{Q}_{(h_1)}) \cap \pi(\tilde{Q}_{(h_2)})
 \subseteq \pi(\tilde{Q}_{(g)})
\]
and that
\[
\pi(\tilde{Q}_{(h_1)}) \cap \pi(\tilde{Q}_{(h_2)}) =
\pi(\tilde{Q}_{(g^\prime)}).
\]
Hence, the intersection of embeddings of two sectors is again a
sector.

\end{proof}

Note that the above lemmas require that the orbifold is cyclic.  It
is easy to construct examples where the above fail in the case of
orbifolds that do not have this property.


\section{The Result}
\label{sec-mainresult}

Now, we return to the question of the existence of a nonvanishing
vector field on a closed almost-complex cyclic orbifold $Q$.  Using
the results of Section \ref{sec-structure}, we will prove Theorem
\ref{thrm-mainresult}.


\begin{proof}[Proof of Theorem \ref{thrm-mainresult}]

Note that as $Q$ is almost-complex, $\tilde{Q}$ inherits
an almost complex structure.  In particular, each of the sectors of $Q$
are even-dimensional and oriented.

\noindent (i) $\Rightarrow$ (ii):

Suppose $X$ is a nonvanishing vector field on $Q$.  Then $X$ is a section of
the orbifold tangent bundle $TQ$, and is hence required to be tangent to the
singular strata of $Q$.  Let $\tilde{X}$ be the induced section
of $\widetilde{TQ}$.  Again, as
$T \tilde{Q} = \widetilde{TQ}$, $\tilde{X}$ is a vector field on $\tilde{Q}$.
However, as $X$ must be tangent to each of the $\pi (\tilde{Q}_{(g)})$ in $Q$,
and as $X$ does not vanish, $\tilde{X}$ is clearly
a nonvanishing vector field on $\tilde{Q}$.

\noindent (ii) $\Rightarrow$ (iii):

Fix a $(g) \in T$ and choose a nonvanishing vector field on
$\tilde{Q}$; we let $X_{(g)}$ denote the restriction of this vector
field to the sector $\tilde{Q}_{(g)}$.  Form the space of sectors
$\widetilde{\tilde{Q}_{(g)}}$ of $\tilde{Q}_{(g)}$ and let
$\widetilde{X_{(g)}}$ denote the induced nonvanishing vector field
on $\widetilde{\tilde{Q}_{(g)}}$.  Clearly, the orbifold index of
$\widetilde{X_{(g)}}$ is zero, and hence by Equation \ref{eq-myph}
(see also Corollary 4.4.4 of \cite{mythesis}), the Euler
characteristic of the underlying space of $\tilde{Q}_{(g)}$ is zero.
As $(g)$ was arbitrary, the Euler characteristic of the underlying
space of each sectors is zero.

\noindent (iii) $\Rightarrow$ (iv):

Suppose all of the Euler characteristics of the sectors of $Q$
vanish. Let $(h)$ be an element that is minimal with respect to the
relation on $T$ so that by Lemma \ref{lemm-minmnflds},
$\tilde{Q}_{(h)}$ is a manifold with the trivial action of a finite
group $C(h) = G_p$, where $G_p$ is the local group with $h \in G_p$.
It is clear from Equation \ref{eq-2ndgb} that
\[
\begin{array}{rcl}
    \chi_{orb}( \tilde{Q}_{(h)} )
        &=& \frac{1}{|G_p|} \chi(\mathbb{X}_{\tilde{Q}_{(h)}} ) \\\\
        &=& 0.
\end{array}
\]
Therefore, the orbifold Euler-Satake characteristics of all of the
minimal elements $\tilde{Q}_{(h)}$ vanish.

Now, fix $(g) \in T$, and suppose that for each $(h) < (g)$, the
orbifold Euler characteristic of $\tilde{Q}_{(h)}$ is zero.  Let
${\mathcal T}$ be a simplicial decomposition for $\tilde{Q}_{(g)}$
such that the isomorphism class of the isotropy group of each point
on the interior of a simplex is constant (see Section
\ref{sec-defs}). Again, for each simplex  $\sigma \in {\mathcal T}$,
we let $m_{\sigma}$ denote the order of the isotropy group of the
points in the interior of $\sigma$.

Now, the Euler-Satake characteristic of $\tilde{Q}_{(g)}$ is given
by
\[
    \chi_{orb}(\tilde{Q}_{(g)}) =   \sum\limits_{\sigma \in
    {\mathcal T}} (-1)^{\mbox{dim}\, \sigma} \frac{1}{m_\sigma}.
\]
Note that the (isomorphism class of the) group generated by $g$ is
contained in the isotropy group of each point in $\tilde{Q}_{(g)}$.
However, any point of $\tilde{Q}_{(g)}$ whose isotropy group is
strictly larger is a singular point of $(\tilde{Q}_{(g)})_{red}$,
and hence is contained in the embedding of a $\tilde{Q}_{(h)}$ for
$(h) < (g)$ by Lemma \ref{lemm-singinsector}. For each simplex
$\sigma$  not completely contained in the image of such an
embedding, the isotropy group is isomorphic to $\langle g \rangle$.
Separating the terms, we have
\begin{equation}
\label{eq-2sum.esc}
\begin{array}{rcl}
    \chi_{orb}(\tilde{Q}_{(g)}) &=&   \sum\limits_{\sigma : m_\sigma = |g|} (-1)^{\mbox{dim}\, \sigma} \frac{1}{m_\sigma}
            + \sum\limits_{\sigma : m_\sigma > |g|} (-1)^{\mbox{dim}\, \sigma}
            \frac{1}{m_\sigma}                  \\\\
                                &=&
            \frac{1}{|g|} \sum\limits_{\sigma : m_\sigma = |g|} (-1)^{\mbox{dim}\, \sigma}
            + \sum\limits_{\sigma : m_\sigma > |g|} (-1)^{\mbox{dim}\, \sigma}
            \frac{1}{m_\sigma}
\end{array}
\end{equation}
The second sum is over all simplices contained in embeddings of
$\tilde{Q}_{(h)}$ for $(h) < (g)$ into $\tilde{Q}_{(g)}$.  While the
images of these embeddings need not be disjoint, the intersection of
any two sectors is again a sector by Lemma \ref{lemm-singinsector}.
Hence, applying Equation \ref{eq-addativeesk}, this sum can be
expressed as a $\Z$-linear combination of Euler-Satake
characteristics of sectors $\tilde{Q}_{(h)}$ with $(h) < (g)$. Using
the fact that each such Euler-Satake characteristic vanishes by the
inductive hypothesis, the second sum is clearly zero.

Similarly, the sum
\begin{equation}
\label{eq-2ndterm.ec}
    \sum\limits_{\sigma : m_\sigma > |g|} (-1)^{\mbox{dim}\, \sigma}
\end{equation}
is the Euler characteristic of the underlying space of the image of
embeddings of sectors.  In the same way, this sum can be expressed
as a $\Z$-linear combination of Euler characteristics of underlying
spaces of sectors $\tilde{Q}_{(h)}$ with $(h) < (g)$.  By
hypothesis, each sector has underlying space with Euler
characteristic 0, and therefore the sum in Equation
\ref{eq-2ndterm.ec} is zero.  Hence,
\begin{equation}
\label{eq-sumequal}
\begin{array}{rcl}
    \sum\limits_{\sigma : m_\sigma > |g|} (-1)^{\mbox{dim}\, \sigma}
\frac{1}{m_\sigma}
        &=&            0                       \\\\
        &=&
            \sum\limits_{\sigma : m_\sigma > |g|} (-1)^{\mbox{dim}\, \sigma}
            \\\\
        &=& \frac{1}{|g|} \sum\limits_{\sigma : m_\sigma > |g|} (-1)^{\mbox{dim}\, \sigma}
\end{array}
\end{equation}

Now, returning to Equation \ref{eq-2sum.esc}, we have
\[
\begin{array}{rcl}
    \chi_{orb}(\tilde{Q}_{(g)})
                                &=&
            \frac{1}{|g|} \sum\limits_{\sigma : m_\sigma = |g|} (-1)^{\mbox{dim}\, \sigma}
            + \sum\limits_{\sigma : m_\sigma > |g|} (-1)^{\mbox{dim}\, \sigma}
            \frac{1}{m_\sigma}              \\\\
                                &=&
            \frac{1}{|g|} \sum\limits_{\sigma : m_\sigma = |g|} (-1)^{\mbox{dim}\, \sigma}
            + \frac{1}{|g|} \sum\limits_{\sigma : m_\sigma > |g|} (-1)^{\mbox{dim}\, \sigma}
                                              \\\\
            &&\mbox{(by Equation \ref{eq-sumequal})}\\\\
            &=&
            \frac{1}{|g|} \sum\limits_{\sigma \in
    {\mathcal T}} (-1)^{\mbox{dim}\, \sigma} \frac{1}{m_\sigma}
                                                \\\\
            &=&
            \frac{1}{|g|} \chi(\mathbb{X}_{\tilde{Q}_{(g)}} )
\end{array}
\]
which is zero by hypothesis. Therefore, by induction, all of the
orbifold Euler-Satake characteristics of the sectors vanish.

\noindent (iv) $\Rightarrow$ (i):

Suppose the Euler-Satake characteristic of each sector of $Q$ is
zero. Adorn $Q$ with a Riemannian metric and extend it in the
natural way to a Riemannian metric on $\tilde{Q}$.  We construct a
nonvanishing vector field $X$ on $Q$ recursively.

Start with the minimal $(h) \in T$.  Each of the corresponding
$\tilde{Q}_{(h)}$ are manifolds with the trivial action of a finite
group by Lemma \ref{lemm-minmnflds}, and by hypothesis, the
Euler-Satake characteristic of each of these manifolds is zero.
Therefore, as the Euler characteristic of each of the underlying
spaces is clearly also zero, it is well know that each
$\tilde{Q}_{(h)}$ admits a nonvanishing vector field. We choose such
a vector field $X_{(h)}$ on each of these minimal sectors.  However,
we require that these vector fields agree on $Q$ in the following
sense: if $\tilde{Q}_{(h_1)} \equiv \tilde{Q}_{(h_2)}$, so that
$\pi(\tilde{Q}_{(h_1)}) = \pi(\tilde{Q}_{(h_2)})$ in $Q$, then the
vector fields on $\tilde{Q}_{(h_1)}$ and $\tilde{Q}_{(h_2)}$
correspond to the same vector field on $Q$.  We can accomplish this
by choosing $X_{(h_1)}$ and defining $X_{(h_2)}$ to be
\[
    X_{(h_2)} :=
        \pi^\ast [(\pi_{(h_1)})^{-1}
            ]^\ast X_{(h_1)}
\]
(see Equation \ref{eq-restrictdiffeo}). In this way, we define $X$
on one representative of each $\equiv$-equivalence class and extend
the definition compatibly to the remaining members of the (finite)
equivalent class.

Now, fix some $(g) \in T$, and suppose that such a nonvanishing vector field
has been given on each $\tilde{Q}_{(h)}$ with $(h) < (g)$ (in such a way that
they agree when $(h_1) \equiv (h_2)$ as described above).  Let
$B:= \{ (p,  (g)) : \exists \: (h) < (g)$ with $(p, (h)) \in  \tilde{Q}_{(h)} \}$,
i.e. the set of all points $(p, (g)) \in \tilde{Q}_{(g)}$ such that
$p$ is fixed by a representative
of some $(h)$ with $(h) < (g)$.  For each such point, $(p, (h))$ is a point in
$\tilde{Q}_{(h)}$ for such an $(h)$, so that $X_{(h)}$ is defined on this
$\tilde{Q}_{(h)}$.  By inverting restrictions of the map $\pi$ and pulling
back the vector field on subsets of $\tilde{Q}_{(g)}$ as above, we define
a vector field on $B$.

Choose a finite set of orbifold charts that cover the (compact) set
$B$. In each chart, we extend the vector field to a parallel vector
field in an open set $W$ containing $B$.  Recall that by Lemma
\ref{lemm-singinsector}, each of the singular points of
$(\tilde{Q}_{(g)})_{red}$ occur as fixed-points of some $h$ with
$(h) < (g)$.  Therefore, $W_{red}$ contains an open neighborhood of
each of the singular points, and $(\tilde{Q}_{(g)})_{red} \backslash
W_{red}$ contains only regular points of $(\tilde{Q}_{(g)})_{red}$.
We extend the vector field to all of $(\tilde{Q}_{(g)})_{red}$ in
such a way that the extension has only isolated zeros and note that
this clearly defines a vector field on $\tilde{Q}_{(g)}$. With this,
we may amalgamate the zeros using well-known techniques (see e.g.
\cite{gp}) by finding a chart with trivial group action that
contains multiple zeros in the image of a compact set.  Such a chart
can be given by choosing a simple path that passes through two zero
points whose image does not intersect $B$ or any of the other zero
points and taking a sufficiently small tubular neighborhood of that
path. Hence, we need not change the vector field on $B$. Moreover,
recall that each of the sectors are even-dimensional, so that the
codimension of the image of each such $\tilde{Q}_{(h)}$ is at least
2; in particular, the preimage of the set $B$ does not separate a
connected uniformized set.

Using this technique, we extend the vector fields $X_{(g)}$ to larger and
larger sectors until we have defined a nonvanishing vector field $X$ on
$\tilde{Q}_{(1)} \cong Q$.

\noindent (i) $\Rightarrow$ (v):

Let $X$ be a nonvanishing vector field on $Q$, and let $\tilde{X}$
be the induced section of $\widetilde{TQ}$.  Restricted to each
sector $\tilde{Q}_{(g)}$, the First Poincar\'{e}-Hopf Theorem for
orbifolds (see \cite{satake2} and \cite{mythesis}) implies that the
integral of the Euler curvature form $E(\Omega)$ defined with
respect to a connection $\omega$ on $\tilde{Q}_{(g)}$ with curvature
$\Omega$ is zero.  Hence, as this closed top form is a
representative of the term in $e_{orb}(Q)$ corresponding to $(g) \in
T$, this term must be zero.  As this is true for each $(g) \in T$,
the orbifold Euler class vanishes.

\noindent (v) $\Rightarrow$ (iv):

Suppose the orbifold Euler class $e_{orb}(Q)$ vanishes.  This implies that
the Euler curvature form $E(\Omega)$ of each of the sectors
$\tilde{Q}_{(g)}$ of $Q$ has integral 0 over the corresponding sector.
By the First Gauss-Bonnet Theorem
for orbifolds (see \cite{satake2} and \cite{mythesis}), the Euler-Satake
characteristics of the sectors all vanish.

\end{proof}

In the case that $Q$ is a manifold, the space of sectors is simply
$Q$ itself, and the Euler characteristic and Euler-Satake
characteristic coincide. Therefore, (i) and (ii) are the same
statement, as are (iii) and (iv).  Additionally, the orbifold Euler
class $e_{orb}(Q)$ reduces to the ordinary Euler class, so that
Theorem \ref{thrm-mainresult} states that $Q$ admits a nonvanishing
vector field if and only if its Euler characteristic vanishes, which
is equivalent to its Euler class vanishing.  Therefore, this theorem
can be viewed as the generalization of the `Hairy Ball Theorem' to
the case of almost-complex orbifolds.

We note that the requirement that $Q$ is almost-complex is crucial,
and not simply required so that the Chen-Ruan cohomology groups are
defined.  In particular, if $Q$ is not almost-complex, then although
all of the singular sets of a sector $\tilde{Q}_{(g)}$ must have
codimension at least 2, it is not necessary that the image of each
$\tilde{Q}_{(h)}$ with $(h) < (g)$ must have codimension at least 2
in $\pi(\tilde{Q}_{(g)})$ (note that the image of the
$\tilde{Q}_{(h)}$ may contain regular points for $\tilde{Q}_{(h)}$
as well as singular points).  Hence, the space formed by removing a
copy of $\tilde{Q}_{(h)}$ from a uniformized set in
$\tilde{Q}_{(g)}$ need not be connected, contributing an additional
obstruction to the amalgamation of zeros of a vector field on
$\tilde{Q}_{(g)}$.

\bibliographystyle{amsplain}

\end{document}